\documentclass[12pt]{amsart}
\usepackage{latexsym,amsmath,amsthm,amssymb,amscd}
\newtheorem{thm}{Theorem}[section]
\newtheorem{lem}[thm]{Lemma}

\newtheorem{remark}[thm]{Remark}
\newtheorem{proposition}[thm]{Proposition}
\newtheorem{example}[thm]{Example}
\newtheorem{defn}[thm]{Definition}

\begin{document}
\def\Ad{{\rm Ad}}
\def\diag{{\rm diag}}
\def\End{{\rm End}}
\def\Fr{{\rm Fr}}
\def\Gal{{\rm Gal}}
\def\SU{{\rm SU}}
\def\SL{{\rm SL}}
\def\GL{{\rm GL}}
\def\Aff{{\rm Aff}}
\def\Id{{\rm I}}
\def\Ind{{\rm Ind}}
\def\Norm{{\rm Norm}}
\def\Nrd{{\rm Nrd}}
\def\Stab{{\rm Stab}}
\def\Id{{\rm Id}}
\def\fa{{\mathfrak a}}
\def\haut{{\rm ht}}
\def\imag{{\rm im}}
\def\integers{{\mathfrak o}}
\def\fR{{\mathfrak R}}
\def\LieN{{\rm n}}
\def\re{{\rm Re}}
\def\Rep{{\rm Rep}}
\def\WDRep{{\rm WDRep}}
\def\rec{{\rm rec}}
\def\Hom{{\rm Hom}}
\def\Ind{{\rm Ind}}
\def\cInd{{\rm c\!\!-\!\!Ind}}
\def\pr{{\rm pr}}
\def\supp{{\rm supp}}
\def\Mod{{\rm Mod}}
\def\ii{{\rm i}}
\def\rr{{\rm r}}
\def\Pbar{{\rm \overline P}}
\def\Nbar{{\rm \overline N}}
\def\sw{{\rm sw}}
\def\Sw{{\rm Sw}}
\def\Ar{{\rm A}}
\def\Sym{{\rm Sym}}
\def\Cusp{{\rm Cusp}}
\def\Irr{{\rm Irr}}
\def\Irrt{{{\rm Irr}^{\rm t}}}
\def\Psit{{{\Psi}^{\rm t}}}
\def\temp{{{\rm t}}}
\def\simple{{{\rm s}}}
\def\Irru{{{\rm Irr}^{\rm u}}}
\def\mes{{\rm mes}}
\def\Tor{{\rm Tor}}
\def\Nor{{\rm N}}
\def\ie{{\it i.e.,\,}}
\def\tr{{\rm tr}}
\def\Omegat{{{\Omega}^{\rm t}}}
\def\cD{{\mathcal D}}
\def\cE{{\mathcal E}}
\def\cH{{\mathcal H}}
\def\cL{{\mathcal L}}
\def\cO{{\mathcal O}}
\def\sA{{\mathcal A}}
\def\sV{{\mathcal V}}
\def\spin{{\rm sp}}
\def\SL{{\rm SL}}
\def\Sp{{\rm Sp}}
\def\SU{{\rm SU}}
\def\red{{\rm nd}}
\def\CC{{\mathbb C}}
\def\QQ{{\mathbb Q}}
\def\RR{{\mathbb R}}
\def\ZZ{{\mathbb Z}}
\def\Afr{{\mathfrak A}}
\def\Bfr{{\mathfrak B}}
\def\Cfr{{\mathfrak C}}
\def\Ofr{{\mathfrak O}}
\def\av{{|\textrm{ }|}}
\def\FF{{F^\times}}
\def\F{{\mathbb F}}
\def\Z{{\mathbb Z}}
\def\ZZ{{\mathbb{Z}^\times}}
\def\R{{\mathbb{R}}}
\def\RR{{\mathbb{R}^\times}}
\def\C{{\mathbb C}}
\def\CC{{\mathbb{C}^\times}}
\def\q{{/\!/}}
\title{Base change and $K$-theory for $\GL(n)$ }
\author{Sergio Mendes and Roger Plymen}
\maketitle

S. Mendes, ISCTE, Av. das For\c{c}as Armadas, 1649-026, Lisbon,
Portugal. Email: sergio.mendes@iscte.pt

R.J. Plymen,  School of Mathematics, Manchester University,
Manchester M13 9PL, England. Email: plymen@manchester.ac.uk

\begin{abstract}
Let $F$ be a nonarchimedean local field and let $G = \GL(n) =
\GL(n,F)$.  Let $E/F$ be a finite Galois extension. We investigate
base change $E/F$ at two levels: at the level of algebraic
varieties, and at the level of $K$-theory.  We put special
emphasis on the representations with Iwahori fixed vectors,  and
the tempered spectrum of $\GL(1)$ and $\GL(2)$. In this context,
the prominent arithmetic invariant is the residue degree $f(E/F)$.
\end{abstract}

Base change. K-theory.


\section{Introduction}

The domain of definition of the classical modular forms (the upper
half plane) is a homogeneous space $\mathbb{H} = \{z \in
\mathbb{C}: \Im z > 0\}$ of the reductive group $G(\mathbb{R}) =
\GL(2,\mathbb{R})$:
\[
\mathbb{H} = \GL(2,\mathbb{R})/O(2,\mathbb{R}) \cdot Z\] where $Z$
is the centre of $G(\mathbb{R})$ and $O(2,\mathbb{R})$ is the
orthogonal group. Each modular form admits a lift $\tilde{f}$ to
the group $\GL(2,\mathbb{R})$ and thence to the adele group
$\GL(2,\mathbb{A})$.

The action of $\GL(2,\mathbb{A})$ on $\tilde{f}$ by right
translation defines a representation $\pi = \pi_f$ of the group
$\GL(2,\mathbb{A})$ in the space of smooth complex-valued
representations on $\GL(2,\mathbb{A})$, for which
\[
\left(\pi(h)\widetilde{f}\right)(g) = \tilde{f}(gh)\] for all $g,h
\in \GL(2,\mathbb{A})$.   If $\pi_f$ is irreducible then one has
an infinite tensor product representation
\[
\pi = \bigotimes_v \pi_v\] where the $\pi_v$ are representations
of the local groups $\GL(2,\mathbb{Q}_v)$ with $v = p$ or
$\infty$.

Let $F$ be a local nonarchimedean field, so that $F$ is either a
finite extension of $\mathbb{Q}_p$ or is a local function field
$\mathbb{F}_q((x))$.  The cardinality of the residue field $k_F$
will be denoted $q_F$.  If $F = \mathbb{Q}_p$ then $q_F = p$. If
$F = \mathbb{F}_q((x))$ then $q_F = q$.

Now let $G = \GL(n) = \GL(n,F)$. Brodzki and Plymen \cite{BP},
working directly with $L$-parameters, equipped the smooth dual of
$\GL(n)$ with a complex structure. In the smooth dual of $\GL(n)$,
especially important are the representations with Iwahori fixed
vectors. This part of the smooth dual has the structure of the
extended quotient $(\CC)^n\q S_n$. This is a smooth complex affine
algebraic variety denoted $\mathfrak{X}_F$.

Let $E/F$ be a finite Galois extension of $F$. We recall that the
domain of an $L$-parameter of $\GL(n,F)$ is the local Langlands
group
\[\mathcal{L}_F: = W_F \times \SL(2,\mathbb{C})\] where $W_F$ is the
Weil group of $F$. Base change is defined by restriction of
$L$-parameter from $\mathcal{L}_F$ to $\mathcal{L}_E$.
 We prove, in section 3, that base change
 \[
 \mathfrak{X}_F \longrightarrow \mathfrak{X}_E\] is a finite morphism of algebraic
 varieties.

An $L$-parameter $\phi$ is \emph{tempered} if $\phi(W_F)$ is
bounded \cite[\S 10.3]{Borel}.   Base change therefore determines
a map of \emph{tempered duals}.  In the rest of this article, we
investigate this map at the level of K-theory.

Let $G(F)= \GL(n,F)$. Let $C^{*}_{r}(G)$ denote the reduced
$C^*$-algebra of $G$. According to the Baum-Connes correspondence,
we have a canonical isomorphism \cite{BHP,La}
$$\mu_{F} : K^{top}_{*}(\beta^1 G(F)) \rightarrow K_{*}C^{*}_{r}(G(F))$$
where $\beta^1 G(F)$ denotes the enlarged building of $G(F)$.

In noncommutative geometry, isomorphisms of $C^*$-algebras are too
restrictive to provide a good notion of isomorphism of
noncommutative spaces, and the correct notion is provided by
strong Morita equivalence of $C^*$-algebras; this point is
emphasized in \cite[p.409]{MP}.   In the present context, the
noncommutative $C^*$-algebra  $C^{*}_{r}(G(F))$ is strongly Morita
equivalent to the commutative $C^*$-algebra $C_{0}(Irr^{t}\,G(F))$
where $Irr^{t}\,G(F)$ denotes the tempered dual of $G(F)$, see
\cite{Pl1}.  As a consequence of this, we have
$$K_{*}C^{*}_{r}(G(F))\cong K^{*}Irr^{t}\,G(F).$$
This leads to the following formulation of the Baum-Connes
correspondence:
$$ K^{top}_{*}(\beta^{1}G(F))\cong K^{*}Irr^{t}\,G(F).$$
This in turn leads to the following diagram

\[\CD
K^{top}_{*}(\beta^{1}G(E)) @> {\mu_{E}}> {}>K^{*}Irr^{t}(G(E))\\
@V{ }VV @VV{b_{E/F}^{*}} V \\ K^{top}_{*}(\beta^{1}G(F)) @>{}>
{\mu_{F}}> K^{*}Irr^{t}(G(F)).
\endCD
\]
where the left-hand vertical map is the unique map which makes the
diagram commutative.

In this paper we focus on the right-hand vertical map. Section 4
contains some partial results.  In sections 5 and 6, we focus on
$\GL(1)$ and $\GL(2)$. We need some crucial results of Bushnell
and  Henniart \cite{BH1}: the exact references are given in
section $6$. In section $6$, the local field $F$ has
characteristic $0$ and $p \neq 2$. The $K$-theory map induced by
unramified base change for \emph{totally ramified cuspidal
representations with unitary central character} is described in
Theorem \ref{K*(BC) GL(2,F)}.

In conformity with the recent book by Bushnell and Henniart
\cite{BH3}, we will consistently write \emph{cuspidal
representation} instead of \emph{supercuspidal representation}.

We would like to thank Guy Henniart for his help on many
occasions, and for his prompt replies to emails.

Sergio Mendes is supported by Funda\c{c}\~{a}o para a Ci\^{e}ncia
e Tecnologia, Terceiro Quadro Comunit\'{a}rio de Apoio,
SFRH/BD/10161/2002.

\section{Base change formula for quasicharacters} Let $F$ be a local nonarchimedean field.
Such a field has an intrinsic norm, denoted $mod_F$ in
\cite[p.4]{We1}. We will write \[ |x|_F = mod_F(x).\] The
valuation $val_F$ is then uniquely determined by the equation
\[ |x|_F = q_F^{- val_F(x)}\]where $q_F$ is the cardinality of the
residue field $k_F = \mathfrak{o}_F/\mathfrak{p}_F$. Here
$\mathfrak{o}_F$ denotes the ring of integers and $\mathfrak{p}_F$
its maximal ideal.

In this section, we review standard material on base change for
quasicharacters. Let $E/F$ be a finite Galois extension, and let
the corresponding Weil groups be denoted $W_E, W_F$.   We have the
standard short exact sequence
\[
1 \to I_E \to W_E \stackrel{d_E}\to \mathbb{Z} \to 0.
\]
Let $Art_E^{-1} : W^{ab}_E \cong E^{\times}$, let $\beta_E : W_E
\to W_E^{ab}$ and let
\[
\alpha_E = Art_E^{-1}\circ \beta_E : W_E \to E^{\times}.
\]

\begin{lem}\label{lemma1}
We have \[N_{E/F}(\alpha_E(w)) = \alpha_F(w)\] for all $w \in W_E
\subset W_F$.\end{lem}
\begin{proof} See  \cite[1.2.2]{Ta1} \end{proof}

\begin{lem}\label{lemma2} We have
\[
f \cdot val_E = val_F \circ N_{E/F}.
\]
\end{lem}
\begin{proof} See \cite[VIII.1, p.139]{We1}.\end{proof}
\begin{lem}\label{lemma3}
We have
\[d_E = - val_E \circ \alpha_E
\]
\end{lem}
\begin{proof}
We have $W_E = \sqcup I_F\Phi_E^n$. Then we have
\begin{eqnarray*} val_E(\alpha_E(x\Phi_E^n)) & = &
val_E(\varpi_E^{-n})\\ & = & -n
\\ & = & -d_E(x\Phi_E^n)\end{eqnarray*} for all $x \in I_E$.
\end{proof}

\begin{lem}\label{lemma4}
Let $w \in W_E \subset W_F$. Then we have
\[
f \cdot d_E(w) = d_F(w).
\]
\end{lem}
\begin{proof} By Lemmas \ref{lemma1}, \ref{lemma2} and \ref{lemma3} we have
\begin{eqnarray*}
d_F(w) & = & - val_F(\alpha_F(w)) \\ & = & -
val_F(N_{E/F}(\alpha_E(w))\\ & = & - f \cdot val_E(\alpha_E(w))\\
& = & f \cdot d_E(w). \end{eqnarray*}
\end{proof}

Now an unramified quasicharacter $\psi$ of $W_E$ is given by the
following simple formula:
\[
\psi(w) = z^{d_E(w)}
\]
where $z \in \mathbb{C}^{\times}$.  The base change formula for a
quasicharacter $\chi$ of $W_F$ is given by
\begin{eqnarray}\label{base change quasicharacters}
b_{E/F}(\chi) = \chi_{|W_E}.
\end{eqnarray}
\begin{lem}
Under base change we have \[ b_{E/F}(\psi)(w) = (z^f)^{d_E(w)}.\]
\end{lem}
for all $w \in W_E.$
\begin{proof}
By Lemma \ref{lemma4} we have
\begin{eqnarray*}
b_{E/F}(\psi)(w) & = & z^{d_F(w)} \\ & = & z^{f \cdot d_E(w)} \\ &
= & (z^f)^{d_E(w)}.
\end{eqnarray*}
\end{proof}

We can remember this result with the (informal) equation
\[z_{E/F} = z^f.\]

\section{Representations with Iwahori fixed vectors}

Let $\Psi(W_{F})$ denote the group of unramified quasicharacters
of $W_{F}$.
Then we have
\[\Psi(W_{F})\cong\CC, \quad \quad \psi \mapsto \psi(\varpi_F).\]

Let $\mathcal{L}_F$ denote the local Langlands group:
\[
\mathcal{L}_F: = W_F \times \SL(2,\mathbb{C}).\]

A \textit{Langlands parameter} (or $L$-parameter) is a continuous homomorphism
$$\phi : \mathcal{L}_F \rightarrow \GL(n,\mathbb{C})$$
($\GL(n,\mathbb{C})$ is given the discrete topology) such that
$\phi(\Phi_F)$ is semisimple, where $\Phi_F$ is a geometric
Frobenius in $W_F$. Two Langlands parameters are equivalent if
they are conjugate under $\GL(n,\mathbb{C})$. The set of
equivalence classes of Langlands parameters is denoted
$\Phi(\GL(n))$.

We will use the local Langlands correspondence for $\GL(n)$
\cite{LRS,He,HT}:
\[
\pi_F : \Phi(\GL(n)) \to \Irr(\GL(n)).\]

Consider first the single $L$-parameter
\[\phi =
1 \otimes \tau(j_1) \oplus \cdots \oplus 1 \otimes
\tau(j_k)\]where $\tau_j$ is the $j$-dimensional complex
representation of $\SL(2,\mathbb{C})$, and $j_1 + \cdots + j_k =
n$.  We define the \emph{orbit} of $\phi$ as follows:
\[\mathcal{O}(\phi) =
\{\psi_1 \otimes \tau(j_1) \oplus \cdots \oplus \psi_k \otimes
\tau(j_k): \psi_r \in \Psi(W_F), 1 \leq r \leq k\}/\sim\]where
$\sim$ denotes the equivalence relation of conjugacy in
$\GL(n,\mathbb{C})$.

 In the local Langlands correspondence,
these $L$-parameters correspond precisely to the irreducible
smooth representations of $\GL(n)$ which admit Iwahori fixed
vectors.

Each partition $j_1 + \cdots + j_k = n$ determines an orbit. The
disjoint union of the orbits, one for each partition of $n$,
creates a complex affine algebraic variety with finitely many
irreducible components. This variety is smooth \cite{BP}.  This
variety admits a simple description as an extended quotient, as we
now proceed to explain.

Let $\Gamma$ be a finite group and $X$ a topological space.
Suppose that $\Gamma$ acts on $X$ as homeomorphisms. Define
$$\widetilde{X}=\{(\gamma , x)\in\Gamma\times X : \gamma x=x\}$$
and
$$g.(\gamma , x)=(g\gamma g^{-1},gx),$$
for all $g\in\Gamma$ and $(\gamma , x)\in\widetilde{X}$. Since $(g\gamma g^{-1})(gx)=g(\gamma x)=gx,$
$\Gamma$ acts on $\widetilde{X}$.
\begin{defn} The extended quotient $X/\!/\Gamma$ associated to the action of
$\Gamma$ on $X$ is the quotient space $\widetilde{X}/\Gamma$.
\end{defn}

If $\gamma\in\Gamma$, let $X^{\gamma}$ denote fixed set
$$X^{\gamma}=\{x\in X : \gamma x=x\}$$
and let $Z_{\gamma}$ denote the centralizer of $\gamma$ in
$\Gamma$. Then the extended quotient is the disjoint union
$$X/\!/\Gamma = \bigsqcup X^{\gamma}/Z_{\gamma}$$
where one $\gamma$ is chosen in each $\Gamma$-conjugacy class.

Let $X= (\mathbb{C}^{\times})^{n}$ be the complex torus of
dimension $n$. The symmetric group $\Gamma =S_{n}$ acts on $X$ by
permuting the coordinates. First, we form the ordinary quotient:
\[
\Sym^n(\CC):= (\CC)^n/S_n\] the $n$-fold symmetric product of
$\CC$.

Next, we form the extended quotient
$(\mathbb{C}^{\times})^{n}/\!/S_{n}$. The conjugacy class of
$\gamma\in S_{n}$ determines a partition of $n$. Let the distinct
parts of the partition be $n_1,\ldots,n_l$ with $n_i$ repeated
$r_i$ times so that
\[
r_1n_1 + \cdots + r_ln_l = n.\]

Let \[ z_j = \psi_j(\varpi_F).\]The map
\[
\psi_1 \otimes \tau(n_1) \oplus \cdots \oplus \psi_{r_1 + \cdots +
r_l} \otimes \tau(n_l) \mapsto (z_1,\ldots, z_{r_1 + \cdots +
r_l})\] determines a bijection
\[\mathcal{O}(\phi) \sim \Sym^{r_{1}}(\CC) \times ... \times
\Sym^{r_{l}}(\CC) = X^{\gamma}/Z(\gamma).\]

With a mild abuse of notation, we will write the $L$-parameter
 \[\phi =
\psi_1 \otimes \tau(j_1) \oplus \cdots \oplus \psi_k \otimes
\tau(j_k)\]
 as
\[
z_1 \cdot \tau(j_1) \oplus \cdots \oplus z_k \cdot \tau(j_k).\]

After base change $E/F$, this $L$-parameter becomes
\[
z_1^f \cdot \tau(j_1) \oplus \cdots \oplus z_k^f \cdot
\tau(j_k).\]

\begin{example}\label{example-GL(4)-extended-quotient}
We illustrate this result for $GL(4)$, by computing the
$L$-parameters, the respective orbits and the extended quotient.
For each item, we list the partition, the $L$-parameter and the
orbit:

\begin{itemize}
\item 4+0, $\phi =1 \otimes \tau(4)$, $\mathcal{O}(\phi) \cong \CC$ \\
\item 3+1,   $\phi =1 \otimes \tau(3)\oplus 1\otimes 1$, $\mathcal{O}(\phi) \cong (\CC)^{2}$ \\
\item 2+2,  $\phi =1\otimes \tau(2)\oplus 1\otimes \tau(2)$, $\mathcal{O}(\phi) \cong Sym^{2}(\CC)$ \\
\item 2+1+1,  $\phi =1 \otimes \tau(2)\oplus 1\otimes 1\oplus
1\otimes 1 $, $\mathcal{O}(\phi) \cong
\CC \times Sym^{2}(\CC)$ \\
\item 1+1+1+1,  $\phi =1 \otimes 1\oplus 1 \otimes 1\oplus 1 \otimes 1\oplus 1 \otimes 1$,
$\mathcal{O}(\phi)\cong Sym^{4}(\CC)$ \\
\end{itemize}

and the extended quotient is
\[(\CC)^4\q S_{4}=\CC \sqcup (\CC)^{2} \sqcup
Sym^{2}(\CC) \sqcup \CC \times Sym^{2}(\CC) \sqcup Sym^{4}(\CC).\]
\end{example}

\begin{thm} Let $\mathfrak{X}_F$ be that part of the smooth
dual of $\GL(n,F)$ comprising all representations which admit
Iwahori fixed vectors.  Then $\mathfrak{X}_F$ is a smooth complex
affine algebraic variety, and in fact has the structure of
extended quotient:
\[\mathfrak{X}_F = (\CC)^n\q S_n.\]
Let $E/F$ be a finite Galois extension.   Then base change
\[\mathfrak{X}_F \longrightarrow \mathfrak{X}_E
\]
is a finite morphism of algebraic varieties.  Explicitly, if $z_1,
\ldots,z_r$ are typical coordinates on $\mathfrak{X}_F$, then base
change is given by
\[
(z_1,\ldots,z_r) \mapsto (z_1^f,\ldots,z_r^f).\]
\end{thm}
\begin{proof} A regular map $f:X\rightarrow{Y}$ of affine varieties is
\emph{finite} if $\mathbb{C}[X]$ is integral over $\mathbb{C}[Y]$,
i.e. if the pullback $f^{\sharp} :
\mathbb{C}[Y]\rightarrow\mathbb{C}[X]$ makes $\mathbb{C}[X]$ a
finitely generated $\mathbb{C}[Y]$-module. The map $\mathfrak{X}_F
\longrightarrow \mathfrak{X}_E$ is regular.

Each irreducible component in the algebraic variety
$\mathfrak{X}_F$ is a product of symmetric products.  Denote a
typical symmetric product by $\mathfrak{S}_F$. The coordinate ring
of each symmetric product $\mathfrak{S}_F$ is of the form
\[ \C[t_1,...,t_r,t_1^{-1},\ldots,t_r^{-1}]^{S_r}
\]the ring of invariant Laurent polynomials.  The pullback is
\[ \C[\mathfrak{S}_E]\longrightarrow\C[\mathfrak{S}_F], t_{i}\mapsto
t_{i}^{f}.
\]

Since $\mathbb{C}[\mathfrak{S}_F]$ is finitely generated as a
$\mathbb{C}[\mathfrak{S}_E]$-module, the base change map is
finite.
\end{proof}

\begin{example} The unramified twists of the Steinberg
representation of $\GL(n)$. \end{example}  These representations
correspond, in the local Langlands correspondence, to the orbit of
the single $L$-parameter $1 \otimes \tau(n)$. This creates an
irreducible curve $\mathfrak{Y}_F$ in the smooth dual of
$\GL(n,F)$, in fact $\mathfrak{Y}_F \cong \CC$. Base change $E/F$
is as follows:
\[
\mathfrak{Y}_F \longrightarrow \mathfrak{Y}_E,\quad \quad z
\mapsto z^f.\]

\begin{example} The spherical Hecke algebra.\end{example}
 Let $K = \GL(n,\mathfrak{o}_F)$, and denote by $\mathcal{H}(G//K)$ the
convolution algebra of all complex-valued, compactly-supported
functions on $G$ such that $f(k_{1}xk_{2})=f(x)$ for all
$k_{1},k_{2}$ in $K$. Then $\mathcal{H}(G//K)$ is called the
\textit{spherical} Hecke algebra.  It is a commutative unital
$\mathbb{C}$-algebra.

 Start with the single $L$-parameter
\[
\phi = 1\otimes 1 \oplus \cdots \oplus 1 \otimes 1\] and let
$\mathcal{O}(\phi)$ denote the orbit of $\phi$. We have
\[
\mathcal{O}(\phi) = \{\psi_1 \otimes 1 \oplus \cdots \oplus \psi_n
\otimes 1\}/\sim\] with $\psi_j$ an unramified quasicharacter of
$W_F$, $1 \leq j \leq n$.  Let $T$ be the standard maximal torus
of $\GL(n)$, and let $^LT$ be the standard maximal torus in the
Langlands dual $^LG$: \[^LT \subset \,^LG =
\GL(n,\mathbb{C}).\]Let $W = S_n$ the symmetric group on $n$
letters. Then we have
\[
\mathcal{O}(\phi) \cong \, \mathbb{C}[^LT/W] = \Sym^n(\CC).\] As a
special case of the above, base change $E/F$ induces the following
finite morphism of algebraic varieties:
\[
\Sym^n(\CC) \longrightarrow \Sym^n(\CC), \quad \quad (z_1,
\ldots,z_n) \mapsto (z_1^f, \ldots,z_n^f).\]


In view of the Satake isomorphism \cite{Ca}
\[\mathcal{H}(G//K) \cong \mathbb{C}[^LT/W]\]
we can interpret base change as an explicit morphism of unital
$\mathbb{C}$-algberas. This recovers (and generalizes) the result
in \cite[p.37]{AC}, for we do not require the extension $E/F$ to
be either unramified or cyclic.

\section{$K$-theory computations}

In this section we compute the $K$-theory map for two examples:
the unitary twists of the Steinberg representation of $\GL(n)$,
and certain connected components in the unitary principal series
of $\GL(n)$.

Let $\mathbb{T}$ denote the circle group
\[\mathbb{T} = \{z \in \mathbb{C}: |z| = 1\}\]
and let $\Psi^t(W_F)$ denote the group of unramified
\emph{unitary} characters of $W_F$. Then we have \[\Psi^t(W_F)
\cong \mathbb{T}, \quad \psi \mapsto \psi(\varpi_F)\]where
$\varpi_F$ is a uniformizer in $F$.

Consider first the single $L$-parameter
\[\phi =
\rho \otimes \tau(j_1) \oplus \cdots \oplus \rho \otimes
\tau(j_k).\] In this formula,  $\rho$ is an irreducible
representation of $W_F$,  $\tau_j$ is the $j$-dimensional complex
representation of $\SL(2,\mathbb{C})$,and $j_1 + \cdots + j_k =
n$. We define the \emph{compact orbit} of $\phi$ as follows:
\[\mathcal{O}^t(\phi) =
\{\bigoplus_{r=1}^{k} \psi_r \otimes \rho \otimes \tau(j_r):
\psi_r \in \Psi^t(W_F), 1 \leq r \leq k\}/\sim\]where, as before,
$\sim$ denotes the equivalence relation of conjugacy in
$\GL(n,\mathbb{C})$.

The Steinberg representation $St_G$ has $L$-parameter $1\otimes
\tau(n)$.

\begin{thm} Let \[
\phi = 1 \otimes \tau(n)\] and let $\mathcal{O}^t(\phi)$ be the
compact orbit of $\phi$. Then we have
\[
BC : \mathbb{T} \to \mathbb{T}, \;\;\; z \mapsto z^f.
\]
$(i)$ This map has degree $f$, and so at the level of the $K$-theory
group $K^1$, $BC$ induces the map
\[
\mathbb{Z} \to \mathbb{Z}, \; \alpha_{1} \mapsto f \cdot \alpha_{1}
\]
of multiplication by the residue degree $f$. Here, $\alpha_{1}$ denotes a generator of $K^{1}(\mathbb{T})\cong\mathbb{Z}$.\\
$(ii)$ At the level of the $K$-theory group $K^0$, $BC$ induces the
identity map
\[
\mathbb{Z} \to \mathbb{Z}, \; \alpha_{0} \mapsto\alpha_{0},
\]
where $\alpha_{0}$ denotes a generator of
$K^{0}(\mathbb{T})\cong\mathbb{Z}$.
\end{thm}
\begin{proof}
$(i)$ Since the map has degree $f$ the result follows. $(ii)$ This
is because $\alpha_0$ is the trivial bundle of rank $1$ over
$\mathbb{T}$.
\end{proof}

Next, we define the $L$-parameter $\phi$ as follows:

$$\phi=\rho \otimes{1}\oplus{...}\oplus\rho\otimes{1}$$
where $\rho$ is a \emph{unitary} character of $W_F$.  The unitary
characters of $W_F$ factor through $F^{\times}$ and we have
\[F^{\times} \cong <\varpi_F>\times \mathcal{U}_F.\] We will take
$\rho$ to be trivial on $<\varpi_F>$, and then regard $\rho$ as a
unitary character of $\mathcal{U}_F$.  The group $\mathcal{U}_F$
admits countably many such characters $\rho$.

 In this case, the
compact orbit is the $n$-fold symmetric product of the circle
$\mathbb{T}$:
$$\mathcal{O}^{t}(\phi)\cong\mathcal{O}^{t}(BC(\phi))\cong\mathbb{T}^{n}/S_{n}.$$

\begin{lem}\label{symmetric power homotopy}
The symmetric product $\mathbb{T}^{n}/S_{n}$ has the homotopy type
of a circle.
\end{lem}
\begin{proof} Send the unordered $n$-tuple $[z]=[z_{1},...,z_{n}]$ to the unique
polynomial with roots $z_{1},...,z_{n}$ and leading coefficient $1$
$$[z_{1},...,z_{n}]\mapsto{z^{n}+a_{n-1}z^{n-1}+...+a_{1}z+a_{0}},\textrm{
}a_{0}\neq{0}.$$ We have then
$$Sym^{n}(\mathbb{C}^{\times})\cong\{z^{n}+a_{n-1}z^{n-1}+...+a_{1}z+a_{0}:a_{0}\neq{0}\}\sim_{h}\mathbb{C}^{\times},$$
since the space of coefficients $a_{n-1},...,a_{1}$ is contractible.
Hence
$$Sym^{n}(\mathbb{T})\sim_{h}\mathbb{T}$$
via the map which sends $[z_{1},...,z_{n}]$ to the product
$z_{1}...z_{n}$.
\end{proof}

We recall the local Langlands correspondence
\[
\pi_F: \Phi(\GL(n)) \to \Irr\, \GL(n).\] Let $t =
diag(x_1,\ldots,x_n)$ be a diagonal element in the standard
maximal torus $T$ of $\GL(n)$. Then \[\sigma: t \mapsto
\pi_F(\rho)(x_1 \cdots x_n)\] is a unitary character of $T$. Let
$\chi$ be an unramified unitary character of $T$, and form the
induced representation \[ \Ind_{TU}^G (\chi \otimes \sigma).\]
This is an irreducible unitary representation of $G$.  When we let
$\chi$ vary over all unramified unitary characters of $T$, we
obtain a subset of the unitary dual of $G$. This subset has the
structure of $n$-fold symmetric product of $\mathbb{T}$.

Since $\mathcal{U}_F$ admits countably many unitary characters,
the unitary dual of $G$ contains countably many subspaces (in the
Fell topology), each with the structure $\Sym^n(\mathbb{T})$. We
are concerned with the effect of base change $E/F$ on each of
these compact spaces.

\begin{thm}\label{K-theory spherical} Let $\mathbb{T}^{n}/S_{n}$
denote one of the compact subspaces of the unitary principal
series of $\GL(n)$ currently under discussion. Then we have
\[
BC : \mathbb{T}^{n}/S_{n} \to \mathbb{T}^{n}/S_{n}, \;\;\;
(z_{1},...,z_{n}) \mapsto (z_{1}^{f},...,z_{n}^{f}).
\]

$(i)$ At the level of the $K$-theory group $K^1$, $BC$ induces the
map
\[
\mathbb{Z} \to \mathbb{Z}, \; \alpha_{1} \mapsto f \cdot \alpha_{1}
\]
of multiplication by $f$, where $f$ is the residue degree and $\alpha_{1}$ denotes a generator of $K^{1}(\mathbb{T})=\mathbb{Z}$.\\
$(ii)$ At the level of the $K$-theory group $K^0$, $BC$ induces
the identity map \[\mathbb{Z} \to \mathbb{Z} , \; \alpha_{0}
\mapsto\alpha_{0},
\]
where $\alpha_{0}$ denotes a generator of
$K^{0}(\mathbb{T})=\mathbb{Z}$.
\end{thm}
\begin{proof}
>From Lemma (\ref{symmetric power homotopy}) we have a commutative
diagram
\[\CD
Sym^{n}(\mathbb{T}) @> {BC}> {}>Sym^{n}(\mathbb{T})\\ @V{h }VV
@VV{h} V \\ \mathbb{T} @>{}> {\widetilde{BC}}> \mathbb{T}
\endCD
\]
Here, $BC(z_{1},...,z_{n})=(z_{1}^{f},...,z_{n}^{f})$, $h$ is the
homotopy map $h([z_{1},...,z_{n}])=z_{1}...z_{n}$ and
$\widetilde{BC}$ is the map $z\mapsto{z^{f}}$. Since \[ (z_1
\cdots z_n)^f = z_1^f \cdots z_n^f\] we have
$K^{j}(BC)=K^{j}(\widetilde{BC})$. But $\widetilde{BC}$ is a map
of degree $f$. Therefore,
$$K^{1}(BC)(\alpha_{1})=f.\alpha_{1} \textrm{ and } K^{0}(BC)(\alpha_{0})=\alpha_{0}$$
where $\alpha_{1}$ is a generator of $K^{1}(\mathbb{T})=\mathbb{Z}$
and $\alpha_{0}$ is a generator of $K^{0}(\mathbb{T})=\mathbb{Z}$.
\end{proof}

\section{Base change and $K$-theory for $\GL(1,F)$}

So far we have considered base change as a map of compact orbits.
Now we want to describe base change as a map of the locally
compact Hausdorff spaces
$$BC:\mathcal{A}^{t}_{1}(F)\rightarrow\mathcal{A}^{t}_{1}(E)$$
where $\mathcal{A}^{t}_{1}(F)$ denotes $Irr^{t}\,\GL(1,F)$. From
now on we will change notation and we denote the tempered dual
$Irr^{t}\,\GL(n,F)$ by $\mathcal{A}^{t}_{n}(F)$.

To study the effect of base change on $K$-theory groups we
explicitly compute the functorial base change map $K^{j}(BC)$. We
will use $K$-theory with compact supports and in particular we will
prove that $BC$ is a proper map.

Let $\chi=\av_{F}^{s}\chi_{0}$ be a character of $F^{\times}$, where
$\chi_{0}$ is the restriction of $\chi$ to
$\mathfrak{o}_{F}^{\times}$. We will write from now on
$\chi=z^{\upsilon_{F}(.)}\chi_{0}$ (since
$|x|^{s}_{F}=q_{F}^{-s\upsilon_{F}(x)}$ this is simply a change of
variables $q_{F}^{-s}\in\mathbb{T}\mapsto z\in\mathbb{T}$). We also
denote the group of units $\mathfrak{o}_{F}^{\times}$ by
$\mathcal{U}_{F}$.

If $\chi_{0}$ is a character of $\mathcal{U}_{F}$ then $\chi_{0}$ is
trivial on some $\mathcal{U}_{F}^{m}$. The least $m$ such that
$\chi_{0}=\chi_{|\mathcal{U}_{F}}$ is trivial on
$\mathcal{U}_{F}^{m}$ is called the \textit{conductor} of $\chi$ and
is denoted $c(\chi)$. Note that $\chi_{0}=\chi_{|\mathcal{U}_{F}}$
can be thought as a character of the finite cyclic group
$\mathcal{U}_{F}/\mathcal{U}_{F}^{c(\chi)}$.

It is well known that the parameters
$(z,c(\chi))\in\mathbb{T}\times\mathbb{N}_{0}$ do not completely
determine the character $\chi$. There is a group isomorphism
$$\mathcal{U}_{F}/\mathcal{U}_{F}^{m}\cong\mathcal{U}_{F}/\mathcal{U}_{F}^{1}\times\mathcal{U}_{F}^{1}/\mathcal{U}_{F}^{2}\times ... \times\mathcal{U}_{F}^{m-1}/\mathcal{U}_{F}^{m}.$$
Now, $\mathcal{U}_{F}/\mathcal{U}_{F}^{1}\cong k_{F}^{\times}$
while $\mathcal{U}_{F}^{i-1}/\mathcal{U}_{F}^{i}\cong k_{F}$ for
$i\geq 1$ \cite[Prop. 5.4]{FV}. Note that $k_{F}^{\times}$ is
interpreted as a multiplicative group while $k_{F}$ is interpreted
as an additive group. Since $k_{F}$ has order $q_{F}$ and
$k_{F}^{\times}$ has order $q_{F}-1$ it follows that
$\mathcal{U}_{F}/\mathcal{U}_{F}^{m}$ is a finite cyclic group of
order $(q_{F}-1)q_{F}^{m-1}$ and
$\widehat{\mathcal{U}_{F}/\mathcal{U}_{F}^{m}}\cong\mathcal{U}_{F}/\mathcal{U}_{F}^{m}$
is also finite with the same order.

We have \cite[Lemma 3.4]{FV}
\[
F^{\times} \cong <\varpi_F> \times \mathcal{U}_F\] where
$\varpi_F$ is a uniformizer if $F$.  It follows that
\[
\mathcal{A}^{t}_{1}(F) \cong \mathbb{T} \times
\widehat{\mathcal{U}_F}.\] We will keep in mind the following
enumeration of the countable set $\widehat{\mathcal{U}_F}$: to
each natural number $n$ we attach the finite set of all characters
$\chi \in \widehat{\mathcal{U}_F}$ for which $c(\chi) = n$.  This
enumeration is not canonical.

\begin{proposition}\cite[Proposition 5, p.143]{We1}
The norm map $N_{E/F}:E^{\times}\rightarrow{F^{\times}}$ determines
an open morphism of $E^{\times}$ onto an open subgroup of
$F^{\times}$.
\end{proposition}

It follows that there exist $m,n\in\mathbb{N}_{0}$ such that
\begin{eqnarray}\label{norm of a compact open}
N_{E/F}(\mathcal{U}_{E}^{n})=\mathcal{U}_{F}^{m}.
\end{eqnarray}

A natural question is how to relate the indexes $n$ and $m$. For unramified extensions we have the following result.
\begin{proposition}\cite[Proposition 1, p.81]{Se1}\label{unramified extension norm map}
Let $E/F$ be a finite, separable, unramified extension. Suppose that $k_{F}$ is finite. Then
$$N_{E/F}(\mathcal{U}_{E}^{n})=\mathcal{U}_{F}^{n}, \textrm{ for all }n\geq{0}.$$
\end{proposition}

Apart from unramified extensions, we will consider the cases when the extension is tamely ramified and totally ramified, since the ramification theory is simpler. The case of wildly ramified extensions is more subtle and will not be considered. We now recall some results about ramification theory. Let $E/F$ be a finite Galois extension and $G=Gal(E/F)$. Put
\begin{eqnarray}\label{def. higher ramification groups}
G_{i}=\{\sigma\in{G}:\sigma{x}-x\in\mathfrak{p}_{E}^{i+1}, \textrm{
for all }x\in\mathfrak{o}_{E}\} \textrm{ , } i\geq{-1}.
\end{eqnarray}
The group $G_{i}$ is called the \textit{i-th ramification group} of the extension $E/F$. Altogether, they form a decreasing sequence of subgroups
\begin{eqnarray}\label{filtration of ramification groups}
G_{-1}=G\supset{G_{0}}\supset{G_{1}}\supset{...}\supset{G_{i}}\supset{G_{i+1}}\supset ...
\end{eqnarray}
Denote by $F_{0}/F$ the maximal unramified subextension of $F$ in $E/F$. Note that $F_{0}$ is the intersection of $E$ with the maximal unramified subextension of $\overline{F}/F$, denoted by $F^{ur}$, where $\overline{F}$ is a fixed algebraic closure of $F$. The subgroup $I_{E/F}=Gal(E/F_{0})$ is called the \textit{inertia subgroup} of $Gal(E/F)$ and we have the identification $G_{0}=I_{E/F}$ \cite[Prop. 1, p.62]{Se1}. It follows that $G/G_{0}=Gal(k_{E}/k_{F})$ and $E/F$ is unramified if and only if $G_{0}=\{1\}$.

The quotient group $G_{0}/G_{1}$ is cyclic and has order prime to
the characteristic residue of $E$ \cite[Cor. 1, p.67]{Se1}.
Therefore, the tame ramification is given by the groups $G_{0}$
and $G_{1}$. In particular, the extension $E/F$ is tamely ramified
if and only if $G_{1}=\{1\}$. These results can be summarize in
the following diagram. \setlength{\unitlength}{3mm}
\begin{equation*}
\begin{picture}(18,12)
\put(0,0){\makebox(0,0){$F$}} \put(0,5){\makebox(0,0){$F_{0}$}}
\put(0,10){\makebox(0,0){$E$}}

\put(5,0){\makebox(0,0){$G$}} \put(5,5){\makebox(0,0){$G_{0}$}}
\put(5,10){\makebox(0,0){$G_{1}$}}
\put(8,10){\makebox(0,0){$=\{1\}$}}

\put(0,1){\line(0,1){3}} \put(0,6){\line(0,1){3}}
\put(1,10){\line(1,0){3}}

\put(5,1){\line(0,1){3}} \put(5,6){\line(0,1){3}}
\put(1,0){\line(1,0){3}} \put(1,5){\line(1,0){3}}

\put(-1,7.5){$\scriptstyle e$}
\put(-1,2.5){$\scriptstyle f$}
\end{picture}
\end{equation*}

To every finite (separable) extension $E/F$ of local fields we associate a real function \cite[\S3, p.73]{Se1}
$$\varphi_{E/F}(u)=\int_{0}^{u}\frac{dt}{(G_{0}:G_{t})}\textrm{ , } \textrm{ for all }u\in [-1,+\infty[.$$
\begin{remark}
This is simply extending the definition of the filtration
(\ref{filtration of ramification groups}) indexed by a finite
discrete parameter to a decreasing filtration $\{G_{t}\}_{t\geq -1}$
with a continuous parameter. If $\textrm{ }i-1<t\leq i$ then we
define $G_{t}=G_{i}$. $\varphi_{E/F}$ is a step function and is a
homeomorphism of the interval $[-1,+\infty [$ into itself.
\end{remark}
The inverse $\psi_{E/F}=\varphi^{-1}_{E/F}$ is called the \textit{Hasse-Herbrand function}. We collect some properties of $\psi_{E/F}$ \cite[Prop. 13, p.73]{Se1}:\\
$(i)$ The function $\psi_{E/F}$ is continuous, piecewise linear, increasing and convex.\\
$(ii)$ $\psi_{E/F}(0)=0$.\\
$(iii)$ If $\nu$ is an integer, then $\psi_{E/F}(\nu)$ is also an integer.
\begin{example}
Suppose $E/F$ is unramified. Then $G_{0}=\{1\}$ and we have
$$\varphi_{E/F}(u)=\int_{0}^{u}\frac{dt}{(G_{0}:G_{t})}=u.$$
Therefore, $\psi_{E/F}(x)=x$.

Now, let $E/F$ be a tame extension. Then $|G_{0}|=e$, $G_{1}=\{1\}$, and we have
$$\varphi_{E/F}(u)=\int_{0}^{u}\frac{dt}{(G_{0}:G_{1})}=\int_{0}^{u}\frac{|G_{1}|}{|G_{0}|}dt=\int_{0}^{u}\frac{dt}{e}=u/e.$$
It follows that $\psi_{E/F}(x)=ex$.
\end{example}
\begin{example}(\cite[p.83]{Se1})
If $E/F$ is cyclic, totally ramified with prime degree $p$, then
\begin{displaymath}
\psi(x)=
 \left\{ \begin{array}{ll}
 x & ,x\leq{t}\\
 t+p(x-t) & ,x\geq{t}
\end{array} \right.
\end{displaymath}
where $t$ is such that $G_{t}\neq\{1\}$ and $\{1\}=G_{t+1}=G_{t+2}=...$ .
\end{example}

\begin{proposition}\cite[Corollary 4, p.93]{Se1}\label{totally ramified extension norm map}
Assume that $E/F$ is a Galois extension, totally ramified. Let $\nu$ be a non-negative number and suppose that $G_{\psi(\nu)}=\{1\}$. Then
$$N_{E/F}(\mathcal{U}_{E}^{\psi(\nu)})=\mathcal{U}_{F}^{\nu}.$$
\end{proposition}
We now deduce a similar result for Galois tamely ramified
extensions.
\begin{proposition}\label{tame extension norm map}
Let $E/F$ be a tamely ramified extension. Then
$$N_{E/F}(\mathcal{U}_{E}^{\psi(\nu)})=\mathcal{U}_{F}^{\nu},$$
for all non-negative real number $\nu$.
\end{proposition}
\begin{proof}
Let $F_{0}/F$ be the maximal unramified subextension of $F$ in $E/F$. We have a tower of fields $F\subset F_{0} \subset E$. Then $E/F_{0}$ is a totally tamely ramified extension. Since $G_{1}=\{1\}$, we also have $G_{\psi(\nu)}=\{1\}$, for all $\nu \geq 0$, where $\psi$ denotes the Hasse-Herbrand function $\psi_{E/F}$.

The conditions of Proposition \ref{totally ramified extension norm map} are satisfied and we have
$$N_{E/F_{0}}(\mathcal{U}_{E}^{\psi(\nu)})=\mathcal{U}_{F_{0}}^{\nu}.$$
Since $F_{0}/F$ is unramified, it follows from Proposition \ref{unramified extension norm map} that
$$N_{F_{0}/F}(\mathcal{U}_{F_{0}}^{\nu})=\mathcal{U}_{F}^{\nu}.$$
Finally, by transitivity of the norm map we have
$$N_{E/F}(\mathcal{U}_{E}^{\psi(\nu)})=N_{F_{0}/F}(N_{E/F_{0}}(\mathcal{U}_{E}^{\psi(\nu)}))=N_{F_{0}/F}(\mathcal{U}_{F_{0}}^{\nu})=\mathcal{U}_{F}^{\nu}.$$
\end{proof}

\textbf{Base change for $\GL(1)$ on the \textit{admissible side}}.

\bigskip

Let $E/F$ be a finite Galois extension and let $W_{E}\hookrightarrow W_{F}$ denote the inclusion of Weil groups. Langlands functoriality predicts the existence a commutative diagram
\[\CD
\mathcal{G}_{1}(F) @> {_{F}\mathcal{L}_{1}}> {}>\mathcal{A}^{t}_{1}(F)\\ @V{Res_{E/F}
}VV @VV{BC} V \\ \mathcal{G}_{1}(E) @>{}> {_{E}\mathcal{L}_{1}}>
\mathcal{A}^{t}_{1}(E)
\endCD
\]
where $\mathcal{G}_{1}(F)$ (resp. $\mathcal{G}_{1}(E)$) is the group
of characters of $W_{F}$ (resp. $W_{E}$) and $BC$ is the base change
map. On the \textit{admissible side} base change is given by
\begin{equation}
\begin{matrix}
\mathcal{A}^{t}_{1}(F) & \longrightarrow & \mathcal{A}^{t}_{1}(E) \cr
\chi_{F} & \longmapsto & \chi_{F}\circ N_{E/F} \cr
\end{matrix}
\end{equation}

\begin{lem}\label{lemma conductor norm}
Let $\chi_{F}$ be a character of $F^{\times}$ with conductor $c(\chi_{F})$ and consider the character $\chi_{E}:=\chi_{F}\circ{N_{E/F}}$ of $E^{\times}$. Suppose we have
$$N_{E/F}(\mathcal{U}_{E}^{n})=\mathcal{U}_{F}^{c(\chi_{F})}.$$
Then, $n$ is the conductor $c(\chi_{E})$ of $\chi_{E}$.
\end{lem}
\begin{proof}
We have
$$\chi_{E}(\mathcal{U}_{E}^{n})=\chi_{F}\circ{}N_{E/F}(\mathcal{U}_{E}^{n})=\chi_{F}(\mathcal{U}_{F}^{c(\chi_{F})})=1.$$
Let $r$ be any integer such that $0<r<n$. Then
\begin{eqnarray*}
\mathcal{U}_{E}^{n} & \subset & \mathcal{U}_{E}^{r} \\
N_{E/F}(\mathcal{U}_{E}^{n}) & \subset & N_{E/F}(\mathcal{U}_{E}^{r}) \\
\mathcal{U}_{F}^{c(\chi_{F})} & \subset & \mathcal{U}_{F}^{s}, \\
\end{eqnarray*}
with $c(\chi_{F})>s$. Then
$$\chi_{F}\circ{}N_{E/F}(\mathcal{U}_{E}^{r})=\chi_{F}(\mathcal{U}_{F}^{s})\neq{1},$$
since $c(\chi_{F})$ is the least integer with this property. Therefore, $n=c(\chi_{E})$.
\end{proof}

We may now describe base change as a map of topological spaces.
The unitary dual $\mathcal{A}^t_1(F)$ (resp. $\mathcal{A}^t_1(E)$)
is a disjoint union of countably many circles, parametrized by
characters $\chi \in \widehat{\mathcal{U}_F}$ (resp. $\eta \in
\widehat{\mathcal{U}_F}$):
\[
\mathcal{A}^t_1(F) \cong \bigsqcup \mathbb{T}_{\chi}, \quad
\mathcal{A}^t_1(E) \cong \bigsqcup \mathbb{T}_{\eta}.\]
We recall that $\chi_{E}=\chi_{F}\circ N_{E/F}$ and $c(\chi_{E})$
is the unique integer such that
$$N_{E/F}(\mathcal{U}_{E}^{c(\chi_{E})})=\mathcal{U}_{F}^{c(\chi_{F})}.$$

\begin{thm}\label{base change GL(1) Hasse-Herbrand}
Let $E/F$ be unramified, tamely ramified or totally ramified (in the last case we also require $E/F$ to be cyclic). Then\\
$(1)$ Base change is a proper map.\\
$(2)$ When we restrict base change to one circle, we get the
following:
\[
BC: \mathbb{T}_{\chi_{F}} \to \mathbb{T}_{\chi_E},\quad z \mapsto
z^f\] with $c(\chi_E) = \psi_{E/F}(c(\chi_{F}))$.
\end{thm}

\begin{proof}
$(1)$ Base change maps each circle into another circle. Let $K$ be
a closed arc in $\mathbb{T}_{\eta}$, and let $\eta = \chi_E$. Then
we may write
$$K=\{e^{i\theta}\in\mathbb{T}_{\eta}: \theta_{0}\leq\theta\leq\theta_{1},\theta\in[0,2\pi]\}.$$
The pre-image of this closed arc is
$$BC^{-1}(K)=\{(e^{i\theta}\in\mathbb{T}_{\chi}: \theta_{0}/f\leq\theta\leq\theta_{1}/f,\theta\in[0,2\pi]\}$$
which is a closed arc in $\mathbb{T}_{\chi}$.  It follows that the
pre-image of a compact set is compact.\\
$(2)$ Follows immediately from Lemma \ref{lemma conductor norm}.
\end{proof}

\smallskip

\textbf{$K$-Theory}

\smallskip

Let $E/F$ be a finite Galois extension. The unitary dual of
$\GL(1)$ is a countable disjoint union of circles and so has the
structure of a locally compact Hausdorff space. The base change
map
\begin{equation}
BC : \bigsqcup \mathbb{T}_{\chi} \longrightarrow \bigsqcup
\mathbb{T}_{\eta}
\end{equation}with $\chi \in \widehat{\mathcal{U}_F},\, \eta \in
\widehat{\mathcal{U}_E}$,  is a proper map.

 Each $K$-group is a countably generated free
abelian group:
$$K^{j}(\mathcal{A}_{1}^{t}(F))\cong\bigoplus \mathbb{Z}_{\chi}, \quad K^{j}
(\mathcal{A}_{1}^{t}(E))\cong\bigoplus \mathbb{Z}_{\eta}$$ with
$\chi \in \widehat{\mathcal{U}_F}, \eta \in
\widehat{\mathcal{U}_E}, j = 0,1$, where $\mathbb{Z}_{\chi}$ and
$\mathbb{Z}_{\eta}$ denote a copy of $\mathbb{Z}$. .

There is a functorial map at the level of $K$-theory groups
\begin{equation}
K^{j}(BC) : \bigoplus \mathbb{Z}_{\eta} \longrightarrow \bigoplus
\mathbb{Z}_{\chi},
\end{equation}

Base change selects among the characters of
$\widehat{\mathcal{U}_E}$ those of the form $\chi_{E}=
\chi_{F}\circ N_{E/F}$, where $\chi_{F}$ is a character of
$\widehat{\mathcal{U}_F}$. \begin{thm}\label{K*(BC) GL(1,F)} When
we restrict $K^1(BC)$ to the direct summand $\mathbb{Z}_{\chi_E}$
we get the following map:
\[
\mathbb{Z}_{\chi_E} \longrightarrow \mathbb{Z}_{\chi_F}, \quad x
\mapsto f \cdot x.\] On the remaining direct summands, $K^1(BC) =
0$. When we restrict $K^0(BC)$ to the direct summand
$\mathbb{Z}_{\chi_E}$ we get the following map:
\[
\mathbb{Z}_{\chi_E} \longrightarrow \mathbb{Z}_{\chi_F}, \quad x
\mapsto x.\] On the remaining direct summands, $K^0(BC) = 0$. In
each case, we have $c(\chi_E) = \psi_{E/F}(c(\chi_F))$.
\end{thm}


\section{Base change and $K$-theory for $\GL(2,F)$}


Through this section, $F$ denotes a nonarchimedean local field
with characteristic $0$ and $p \neq 2$.

Let $\mathcal{G}^{0}_{2}(F)$ be the set of equivalence classes of
irreducible $2$-dimensional smooth (complex) representations of
$W_{F}$. Let $\mathcal{A}^{0}_{2}(F)$ be the subset of
$\mathcal{A}^{t}_{2}(F)$ consisting of equivalence classes of
irreducible cuspidal representations of $\GL(2,F)$. The local
Langlands correspondence gives a bijection
$$_{F}\mathcal{L} : \mathcal{G}^{0}_{2}(F) \rightarrow \mathcal{A}^{0}_{2}(F).$$

We recall the concept of \emph{admissible pair} \cite[p.124]{BH3}.
\begin{defn}
Let $E/F$ be a quadratic extension and let $\xi$ be a
quasicharacter of $E^{\times}$. The pair $(E/F,\xi)$ is called
admissible if\\
$(1)$ $\xi$ does not factor through the norm map $N_{E/F}: E^{\times} \to F^{\times}$ and, \\
$(2)$ if $\xi|\mathcal{U}_{E}^{1}$ does factor through $N_{E/F}$,
then $E/F$ is unramified.
\end{defn}
Denote the set of $F$-isomorphism classes of admissible pairs
$(E/F,\xi)$ by $\mathcal{P}_{2}(F)$. According to
\cite[p.215]{BH3}, the map
\begin{equation}
\begin{matrix}
\mathcal{P}_{2}(F) & \longrightarrow & \mathcal{G}^{0}_{2}(F) \cr
(E/F,\xi) &\longmapsto & \Ind_{E/F} \xi \cr
\end{matrix}
\end{equation}
is a canonical bijection, where we see $\xi$ as a quasicharacter
of $W_{E}$ via the class field theory isomorphism $W_{E}\cong
E^{\times}$ and $\Ind_{E/F}$ is the functor of induction from
representations of $W_{E}$ to representations of $W_{F}$.


The tempered dual of $\GL(2)$ comprises the cuspidal
representations with unitary central character, the unitary twists
of the Steinberg representation, and the unitary principal series.

We will restrict ourselves to admissible pairs $(E/F,\xi)$ for
which $\xi$ is a unitary character. This ensures that $\rho: =
\Ind_{E/F}\xi$ is unitary. Therefore $\det(\rho)$ is unitary and
$\mathcal{L}(\rho)$ has unitary central character.

The cuspidal representations of $\GL(2)$ with unitary central
character arrange themselves in the tempered dual as a countable
union of circles.  For each circle $\mathbb{T}$, we select an
admissible pair $(E/F,\xi)$ for which
\[\mathcal{L}(\rho) \in \mathbb{T}\] and label this circle as
$\mathbb{T}_{(E/F,\xi)}$.

We \emph{further restrict ourselves} to admissible pairs
$(E/F,\xi)$ for which $E/F$ is totally ramified.

\begin{thm}\label{base change GL(2) Hasse-Herbrand}
Let $L/F$ be an unramified extension of odd
degree. Then\\
$(1)$ Base change is a proper map,\\
$(2)$ When we restrict base change to one circle we get the
following:
\[
BC: \mathbb{T}_{(E/F,\xi)} \to \mathbb{T}_{(EL/L,\xi_L)},\quad z
\mapsto z^{f(L/F)}\] with $\xi_L =\xi\circ N_{EL/E}$ and $c(\xi_L) =
\psi_{EL/E}(c(\xi))=c(\xi)$.
\end{thm}
\begin{proof}
The proof of $(1)$ is analogous to the proof of Theorem \ref{base
change GL(1) Hasse-Herbrand}.

Each representation $\rho\in\mathcal{G}^0_2(F)$ has a
\emph{torsion number}: the order of the cyclic group of all the
unramified characters $\chi$ for which $\chi\rho\cong\rho$. The
torsion number of $\rho$ will be denoted $t(\rho)$.

Set $\sigma = \Ind_{E/F}\xi, \pi = \mathcal{L}(\sigma)$ and
$\sigma_L = \Ind_{EL/L}\xi_L = \sigma|W_L$.  Then $\sigma$ is
totally ramified, in the sense that $t(\sigma) = 1$, as in the
proof of Theorem 3.3 in \cite[p. 697]{BH1}. The pair
$(EL/L,\xi_L)$ is admissible \cite[Theorem 4.6]{BH1}. We now quote
\cite[Proposition 3.2(7)]{BH1} to infer that
\[
\mathcal{L}(\sigma_L) = b_{L/F}\,\pi.\]

If $L/F$ is unramified then, by \cite[Proposition 7.2, p.153]{Ne},
$EL/E$ is unramified.  For the ramification indices, we have:
$$e_{EL/F}=e_{EL/L}\times e_{L/F}=e_{EL/E}\times e_{E/F}.$$
Since $L/F$ and $EL/E$ are both unramified it follows that
\[e_{EL/L}=e_{E/F} = 2.\]

Since $EL/L$ is a quadratic extension, $EL/L$ is totally ramified.
Therefore $\sigma_L$ is totally ramified, i.e. $t(\sigma_L) = 1$.

Hence, base change maps each circle into another circle and the
map is given by $z \mapsto z^{f(L/F)}$.

Finally, since $EL/E$ is unramified, we have $\psi_{EL/E}(x)=x$
and the result follows.

\end{proof}

Let $L/F$ be a finite unramified Galois extension. The cuspidal
part of the tempered dual of $\GL(2)$ is a countable disjoint
union of circles and so has the structure of a locally compact
Hausdorff space. The base change map
\begin{equation}
BC : \bigsqcup \mathbb{T}_{(E/F,\xi)} \longrightarrow \bigsqcup
\mathbb{T}_{(EL/L,\eta)}
\end{equation}
with $(E/F,\xi)$ an admissible pair, $E/F$ totally ramified and
$\xi$ unitary is a proper map.

 Each $K$-group is a countably generated free
abelian group:
$$K^{j}( \bigsqcup \mathbb{T}_{(E/F,\xi)})\cong\bigoplus \mathbb{Z}_{(E/F,\xi)}, \quad K^{j}
(\bigsqcup \mathbb{T}_{(EL/L,\eta)})\cong\bigoplus
\mathbb{Z}_{(EL/L,\eta)}$$ where $\mathbb{Z}_{(E/F,\xi)}$ and
$\mathbb{Z}_{(EL/L,\eta)}$ denote a copy of $\mathbb{Z}$, $j = 0,1$.

In complete analogy with $\GL(1)$ there is a functorial map at the
level of $K$-theory groups
\begin{equation}
K^{j}(BC) : \bigoplus \mathbb{Z}_{(EL/L,\eta)} \longrightarrow
\bigoplus \mathbb{Z}_{(E/F,\xi)}.
\end{equation}

Base change selects among the admissible pairs $(EL/L,\eta)$ those
of the form $(EL/L,\xi_L)$, where $\xi_{L}=\xi\circ N_{EL/E}$.
\begin{thm}\label{K*(BC) GL(2,F)} When we restrict $K^1(BC)$ to the
direct summand $\mathbb{Z}_{(EL/L,\xi_L)}$ we get the following map:
\[
\mathbb{Z}_{(EL/L,\xi_L)} \longrightarrow \mathbb{Z}_{(E/F,\xi)},
\quad x \mapsto f(L/F) \cdot x.\] On the remaining direct
summands, $K^1(BC) = 0$. When we restrict $K^0(BC)$ to the direct
summand $\mathbb{Z}_{(EL/L,\xi_L)}$ we get the following map:
\[
\mathbb{Z}_{(EL/L,\xi_L)} \longrightarrow \mathbb{Z}_{(E/F,\xi)},
\quad x \mapsto x.\] On the remaining direct summands, $K^0(BC) =
0$.
\end{thm}

\end{document}